\documentclass[a4paper,12pt]{amsart}
\usepackage{amsfonts}
\usepackage{amssymb}
\usepackage{ifthen}
\usepackage{graphicx}
\usepackage{color}
\nonstopmode \numberwithin{equation}{section}
\newtheorem{thm}{Theorem}
\newtheorem{lem}{Lemma}
\newtheorem{cor}{Corollary}[section]

\newtheorem{cl}{Claim}
\newtheorem{ca}{Case}
\newtheorem{sca}{Subcase}
\newtheorem{scl}{Subclaim}
\newtheorem{conj}[equation]{Conjecture}

\theoremstyle{definition}
\newtheorem{defn}{Definition}

\newtheorem{op}[equation]{Open Problem}
\newtheorem{ques}[equation]{Question}
\newtheorem{rem}{Remark}[section]
\newtheorem{exam}[equation]{Example}

\newcounter {own}
\def\theown {\thesection       .\arabic{own}}

\newenvironment{pf}[1][]{%
 \vskip 3mm
 \noindent
 \ifthenelse{\equal{#1}{}}%
  {{\slshape Proof. }}%
  {{\slshape #1.} }%
 }%
{\qed\bigskip}

\newcounter{alphabet}
\newcounter{tmp}
\newenvironment{Thm}[1][]{\refstepcounter{alphabet}%
\bigskip%
\noindent%
{\bf Theorem \Alph{alphabet}}%
\ifthenelse{\equal{#1}{}}{}{ (#1)}%
{\bf .} \itshape}{\vskip 8pt}

\makeatletter
\newcommand{\Ref}[1]{\@ifundefined{r@#1}{}{\setcounter{tmp}{\ref{#1}}\Alph{tmp}}}
\makeatother

\newenvironment{Lem}[1][]{\refstepcounter{alphabet}%
\bigskip%
\noindent%
{\bf Lemma \Alph{alphabet}}%
{\bf .} \itshape}{\vskip 8pt}

\newcommand{\ID}{{\mathbb D}}

\newcommand{\arctanh}{{\operatorname{arctanh}}}

\def\be{\begin{equation}}
\def\ee{\end{equation}}

\newcommand{\bee}{\begin{enumerate}}
\newcommand{\eee}{\end{enumerate}}

\newcommand{\blem}{\begin{lem}}
\newcommand{\elem}{\end{lem}}
\newcommand{\bthm}{\begin{thm}}
\newcommand{\ethm}{\end{thm}}
\newcommand{\bcor}{\begin{cor}}
\newcommand{\ecor}{\end{cor}}
\newcommand{\beg}{\begin{exam}}
\newcommand{\eeg}{\end{exam}}
\newcommand{\begs}{\begin{examples}}
\newcommand{\eegs}{\end{examples}}
\newcommand{\bdefe}{\begin{defn}}
\newcommand{\edefe}{\end{defn}}
\newcommand{\bprob}{\begin{prob}}
\newcommand{\eprob}{\end{prob}}
\newcommand{\bques}{\begin{ques}}
\newcommand{\eques}{\end{ques}}
\newcommand{\bei}{\begin{itemize}}
\newcommand{\eei}{\end{itemize}}
\newcommand{\bcon}{\begin{conj}}
\newcommand{\econ}{\end{conj}}
\newcommand{\bop}{\begin{op}}
\newcommand{\eop}{\end{op}}

\newcommand{\bca}{\begin{ca}}
\newcommand{\eca}{\end{ca}}
\newcommand{\bsca}{\begin{sca}}
\newcommand{\esca}{\end{sca}}

\newcommand{\bcl}{\begin{cl}}
\newcommand{\ecl}{\end{cl}}

\newcommand{\bscl}{\begin{scl}}
\newcommand{\escl}{\end{scl}}

\newcommand{\bcons}{\begin{conjs}}
\newcommand{\econs}{\end{conjs}}
\newcommand{\bprop}{\begin{propo}}
\newcommand{\eprop}{\end{propo}}
\newcommand{\br}{\begin{rem}}
\newcommand{\er}{\end{rem}}
\newcommand{\brs}{\begin{rems}}
\newcommand{\ers}{\end{rems}}
\newcommand{\bo}{\begin{obser}}
\newcommand{\eo}{\end{obser}}
\newcommand{\bos}{\begin{obsers}}
\newcommand{\eos}{\end{obsers}}
\newcommand{\bpf}{\begin{pf}}
\newcommand{\epf}{\end{pf}}
\newcommand{\ba}{\begin{array}}
\newcommand{\ea}{\end{array}}
\newcommand{\beq}{\begin{eqnarray}}
\newcommand{\beqq}{\begin{eqnarray*}}
\newcommand{\eeq}{\end{eqnarray}}
\newcommand{\eeqq}{\end{eqnarray*}}

\newcommand{\ra}{\rightarrow}

\newcommand{\ds}{\displaystyle}

\newcounter{minutes}
\divide\time by 60
\newcounter{hours}
\multiply\time by 60 \addtocounter{minutes}{-\time}

\begin{document}

\bibliographystyle{amsplain}
\title [] {Coefficient estimates, Landau's theorem and Lipschitz-type spaces on planar harmonic mappings}

\def\thefootnote{}
\footnotetext{ \texttt{\tiny File:~\jobname .tex,
          printed: \number\day-\number\month-\number\year,
          \thehours.\ifnum\theminutes<10{0}\fi\theminutes}
} \makeatletter\def\thefootnote{\@arabic\c@footnote}\makeatother

\author{Sh. Chen}
\address{Sh. Chen, Department of Mathematics and Computational
Science, Hengyang Normal University, Hengyang, Hunan 421008,
People's Republic of China.} \email{shlchen1982@yahoo.com.cn}

\author{S.  Ponnusamy $^\dagger $
}
\address{S. Ponnusamy,
Indian Statistical Institute (ISI), Chennai Centre, SETS (Society
for Electronic Transactions and security), MGR Knowledge City, CIT
Campus, Taramani, Chennai 600 113, India. }
\email{samy@isichennai.res.in, samy@iitm.ac.in}

\author{ A. Rasila }
\address{A. Rasila, Department of Mathematics and Systems Analysis, Aalto University, P. O. Box 11100, FI-00076 Aalto,
 Finland.} \email{antti.rasila@iki.fi}

\keywords{Harmonic mappings, harmonic univalent mappings,  Lipschitz space,  harmonic Bloch space, Green's theorem.\\
$
^\dagger$ {\tt Corresponding author. This author is on leave from the Department of Mathematics,
Indian Institute of Technology Madras, Chennai-600 036, India}
}

\begin{abstract}
In this paper, we investigate the properties of locally univalent
and multivalent planar harmonic mappings. First, we discuss the
coefficient estimates and Landau's Theorem for some classes of
locally univalent harmonic mappings, and then we study some
Lipschitz-type spaces for locally univalent and multivalent harmonic
mappings.
\end{abstract}


\maketitle \pagestyle{myheadings} \markboth{ SH. Chen, S. Ponnusamy
and A. Rasila}{Coefficient estimates, Landau's theorem and
Lipschitz-type spaces on planar harmonic mappings }

\section{Introduction }\label{csw-sec1}
Let $D$ be a simply connected subdomain of the complex plane
$\mathbb{C}$. A complex-valued function $f$ defined in $D$ is called
a {\it harmonic mapping}  if and only if both the real and the
imaginary parts of $f$ are real harmonic in $D$.
It is known that every harmonic mapping $f$ defined in $D$ admits a
decomposition $f=h+\overline{g}$, where $h$ and $g$ are analytic.
Since the Jacobian $J_f$ of $f$ is given by
$$J_f=|f_{z}|^2-|f_{\overline{z}}|^2:=|h'|^2-|g'|^2,
$$
$f$ is locally univalent and sense-preserving in  $D$ if
and only if $|g'(z)|<|h'(z)|$ in $D$; or equivalently if
$h'(z)\neq0$ and the dilatation $w=g'/h'$ has the property that
$|w(z)|<1$ in $D$ (see \cite{Lewy}). We refer to
\cite{Clunie-Small-84,Du} for the theory of planar harmonic
mappings.

For $a\in\mathbb{C}$, let $\ID(a,r)=\{z:\, |z-a|<r\}$. In
particular, we use $\mathbb{D}_r$ to denote the disk
$\mathbb{D}(0,r)$ and  $\mathbb{D}$ the open unit disk $\ID_1$. For
harmonic mappings $f$ defined in $\mathbb{D}$, we use the following
standard notations:
$$\Lambda_{f}(z)=\max_{0\leq \theta\leq 2\pi}|f_{z}(z)+e^{-2i\theta}f_{\overline{z}}(z)|
=|f_{z}(z)|+|f_{\overline{z}}(z)|
$$
and
$$\lambda_{f}(z)=\min_{0\leq \theta\leq 2\pi}|f_{z}(z)+e^{-2i\theta}f_{\overline{z}}(z)|
=\big | \, |f_{z}(z)|-|f_{\overline{z}}(z)|\, \big |.
$$
Thus, for a sense-preserving harmonic mapping $f$, one has
$J_f(z)=\Lambda_{f}(z)\lambda_{f}(z).
$


\section{Main results}\label{csw-sec2}

For a harmonic mapping $f$ in $\mathbb{D}$ and $r\in[0,1)$, the
{\it harmonic area function} $S_{f}(r)$ of $f$, counting
multiplicity, is defined by
$$S_{f}(r)=\int_{\mathbb{D}_{r}}J_f(z)\,dA(z),
$$
where $dA$ denotes the normalized area measure on $\mathbb{D}$. In
\cite{KK}, the authors discussed some properties of  harmonic area
functions and proved that a harmonic self-homeomorphism of a disk
does not increase the area of any concentric disk. In this paper, we
discuss  coefficients estimates and the  Landau Theorem for
sense-preserving harmonic mappings having  finite area.
 Let ${\mathcal H}$ denote   the set of all
sense-preserving harmonic mappings $f=h+\overline{g}$ in
$\mathbb{D}$ satisfying the normalization
$h(0)=g(0)=f_{\overline{z}}(0)=0$, where $h$ and $g$ are analytic in
$\mathbb{D}$. We denote ${\mathcal H}(C)$ the class of all mappings
$f=h+\overline{g}\in{\mathcal H}$ with the finiteness condition
$$\ds C:=\sup_{0<r<1} S_{f}(r)<\infty,
$$ where $h$ and $g$ are analytic in $\mathbb{D}$. For any fixed  $\alpha\geq0$,
let ${\mathcal H}_{\alpha}(C)$ denote all  mappings $f\in {\mathcal H}(C)$ with
$f_{z}(0)=\alpha$.

 In order to state our main results, we first consider the class
 ${\mathcal H}(C)$. Throughout the discussion, we assume that $h$ and $g$ have the form
\be\label{CRP-eq-3}
 h(z)=\sum_{n=1}^{\infty}a_{n}z^{n}
~\mbox{ and }~g(z)=\sum_{n=2}^{\infty}b_{n}z^{n}, ~z\in \ID,
\ee
whenever $f=h+\overline{g}$ belongs to ${\mathcal H}(C).$ Our
results are as follows.

\begin{thm}\label{thm-1}
Let  $f\in {\mathcal H}(C)$,
$r_{0}=(\sqrt{5}-1)/2\approx0.618$, and
$$Q(r_{0})=\left\{\frac{(1+r_{0})}{r_{0}^{2}(1-r_{0})}C\right\}^{1/2} \approx3.330\sqrt{C}.
$$
Then   

$$
\begin{cases}
\displaystyle  |a_{1}|\leq\sqrt{2C}
& \mbox{ if } n=1,\\[3mm]
\displaystyle |a_{n}|+|b_{n}|\leq\frac{4Q(r_{0})}{ \pi
r_{0}^{n-1}}\left(1+\frac{1}{n-1}\right)^{n-1}<\frac{4Q(r_{0}){\rm
e}}{ \pi r_{0}^{n-1}} &\mbox{ if }\, n\geq2,
\end{cases}
$$
where
\be\label{defn-e}
{\rm e} =\lim_{n\rightarrow\infty}\left(1+\frac{1}{n}\right)^{n}.
\ee
In the special case that $C=1/2$,  the estimate of $|a_{1}|$ is
sharp and the extreme function is $f(z)=z+\overline{z}^{2}/2$.

\end{thm}

Let $f$ be a sense-preserving harmonic mapping from $\mathbb{D}$
into $\mathbb{C}$. We say that $f$ is a {\it $K$-quasiregular
harmonic mapping} if and only if
$$\frac{\Lambda_{f}(z)}{\lambda_{f}(z)}\leq K,  ~\mbox{ i.e., }~\frac{|f_{\overline{z}}(z)|}{|f_{z}(z)|}\leq \frac{K-1}{K+1},
~\mbox{ for }~z\in \ID,
$$
where $K\geq1$. Moreover, if $f$ is a univalent and $K$-quasiregular
harmonic mapping, then $f$ is called a $K$-quasiconformal harmonic
mapping.

 A harmonic mapping $f$ is called a {\it harmonic Bloch mapping} if and only if
$$\sup_{ z,w\in\mathbb{D},\ z\neq w}\frac{|f(z)-f(w)|}{\sigma(z,w)}< \infty,
$$
where
$$\sigma(z,w)=\frac{1}{2}\log\left(\frac{|{1-\overline{z}w}|+\left|z-w\right|}
{|{1-\overline{z}w}|-\left|z-w\right|}\right)=\arctanh \left
|\frac{z-w}{1-\overline{z}w}\right|
$$
denotes the hyperbolic distance between $z$ and $w$ in $\mathbb{D}$.

In \cite{Co}, Colonna proved that \be\label{eq1} \sup_{
z,w\in\mathbb{D},\ z\neq w}\frac{|f(z)-f(w)|}{\sigma(z,w)}
=\sup_{z\in\mathbb{D}}\{(1-|z|^{2})\Lambda_{f}(z)\}. \ee Moreover,
the set of all harmonic Bloch mappings, denoted by the symbol
$\mathcal{HB}$, forms a complex Banach space  with the norm
$\|\cdot\|$ given by
$$
\|f\|_{\mathcal{HB}}=|f(0)|+\sup_{z\in\mathbb{D}}\{(1-|z|^{2})\Lambda_{f}(z)\}.
$$

For $K$-quasiregular harmonic mappings, we have
\begin{thm}\label{thm-1.x}
Let $f=h+\overline{g}$ be a $K$-quasiregular harmonic mapping in
$\mathbb{D}$ satisfying  $C=\sup_{0<r<1} S_{f}(r)<\infty$, where
$h(z)=\sum_{n=1}^{\infty}a_{n}z^{n}~\mbox{and}~g(z)=\sum_{n=1}^{\infty}b_{n}z^{n}.$
  Then $f\in\mathcal{HB}$ and
$$|a_{n}|+|b_{n}|\leq
\begin{cases}
\displaystyle  \sqrt{CK}
& \mbox{ if } n=1,\\[3mm]
\displaystyle
\frac{4\sqrt{CK}}{\pi}\left(1+\frac{1}{n-1}\right)^{n-1} &\mbox{ if
}\, n\geq2.
\end{cases}
$$
In the special case that $C=K=1$,  the estimate of $|a_{1}|$ is
sharp and the extreme function is $f(z)=z$.

\end{thm}

The following result is obtained as an application of Theorem \ref{thm-1}.

\begin{thm}\label{thm-2}
Let $f\in {\mathcal H}_{\alpha}(C)$
with $0<\alpha< Q(r_{0})$, where $r_{0}$ and $Q(r_{0})$ are the same
as in Theorem {\rm \ref{thm-1}}. Then for $n\geq2$,
$$|a_{n}|+|b_{n}|\leq\frac{1}{nr_{0}^{n-1}Q(r_{0})}\inf_{0<t<1}\left\{\frac{Q(r_{0})^{2}-\alpha^{2}(1-t)^{2}}{t^{n-1}(1-t)}\right\}<
\frac{Q(r_{0})}{ r_{0}^{n-1}}\left(1+\frac{1}{n-1}\right)^{n-1}.
$$
\end{thm}

The classical theorem of Landau asserts the existence of an
universal constant $\rho$ such that every analytic function $f:\,\ID
\ra \ID $ with $f(0)=f'(0)-1=0$ is univalent in the disk
$\mathbb{D}_\rho $ and in addition, the range $f(\mathbb{D}_\rho )$
contains a disk of radius $\rho ^2$. Recently, many authors
considered Landau's theorem for planar harmonic mappings (see for
example,
\cite{HG,HG1,CPW0,CPW1,CPW4,CPW6,CPW7,DN-04,Armen-06,Liu1,Hu}).
Applying Theorem \ref{thm-2}, we obtain the following result, and since
a bounded harmonic mapping in $\mathbb{D}$ has a finite area, we see
that this result is indeed a generalization of \cite[Theorem
3]{HG}.

\begin{thm}\label{thm2.0}
Let $f\in {\mathcal H}_{\alpha}(C)$
with $0<\alpha< Q(r_{0})$, where $r_{0}$ and $Q(r_{0})$ are the same
as in Theorem {\rm \ref{thm-1}}. Define
$$\rho =1-\frac{1}{\sqrt{1+ \alpha /({\rm e} Q(r_{0}) )}} ~\mbox{ and }
~R_{0}=r_{0}\rho\left(\alpha-\frac{{\rm e} Q(r_{0})\rho}{1-\rho}\right).
$$
Then $f$ is univalent in $\mathbb{D}_{\rho}$. Moreover, $f(r_{0}\mathbb{D}_{\rho})$
contains a univalent disk $\mathbb{D}_{\mbox{R}_{0}}.$
\end{thm}

A continuous increasing function $\omega:\, [0,\infty)\rightarrow
[0,\infty)$ with $\omega(0)=0$ is called a {\it majorant} if
$\omega(t)/t$ is non-increasing for $t>0$ (see \cite{D}). Given a
subset $\Omega$ of $\mathbb{C}$, a function $f:\, \Omega\rightarrow
\mathbb{C}$ is said to belong to the {\it Lipschitz space
$L_{\omega}(\Omega)$} if there is a positive constant $M$ such that
\be\label{eq1x}
|f(z)-f(w)|\leq M\omega(|z-w|) ~\mbox{ for all $z,\ w\in\Omega$.}
\ee

For $\delta_{0}>0$ and $0<\delta<\delta_{0}$, we consider the
following conditions on a majorant $\omega$:
\be\label{eq2x}
\int_{0}^{\delta}\frac{\omega(t)}{t}\,dt\leq M \omega(\delta)
\ee
and
\be\label{eq3x}
\delta\int_{\delta}^{+\infty}\frac{\omega(t)}{t^{2}}\,dt\leq M
\omega(\delta),
\ee
where $M$ denotes a positive constant.

A majorant $\omega$ is said to be {\it regular} if it satisfies the
conditions (\ref{eq2x}) and (\ref{eq3x}) (see \cite{D}).

Dyakonov \cite{D} characterized the holomorphic functions of class
$L_{\omega}$ in terms of their modulus. Later in \cite[Theorems
A~\mbox{and}~B]{P}, Pavlovi\'{c} came up with a relatively simple
proof of the results of Dyakonov. Recently, many authors considered
this topic and generalized Dyakonov's results to holomorphic
functions and  harmonic functions of one variable and several
variables  (see
\cite{MVM1,CPW4,D,D1,M3,MM,P, Pav1,Pav2,Pav3}). In this paper, we
first extend \cite[Theorems A and B]{P} to planar harmonic mappings
as follows.

\begin{thm}\label{thm-7}
Let $\omega$ be a  majorant satisfying \eqref{eq2x} and
 $f\in{\mathcal H}$. Then for all
$r\in(0,1)$,
$$f\in L_{\omega}(\mathbb{D}_{r})\Longleftrightarrow
|f|\in L_{\omega}(\mathbb{D}_{r}) \Longleftrightarrow |f|\in
L_{\omega}(\mathbb{D}_{r},\partial \mathbb{D}_{r}),
$$
where $L_{\omega}(\mathbb{D}_{r},\partial \mathbb{D}_{r})$ denotes
the class of continuous functions $F$ on $\mathbb{D}_{r}\cup\partial
\mathbb{D}_{r}$ which satisfy the condition \eqref{eq1x} with some positive
constant $M$, whenever $z\in \mathbb{D}_{r}$ and $w\in\partial
\mathbb{D}_{r}$.
\end{thm}

In \cite{KO}, Korenblum proved the following result.

\begin{Thm}\label{Thm-X}
Let $u$ be a real harmonic function in $\mathbb{D}$ and $f_{r}(\theta)=u(re^{i\theta})$. Then
$$\|f_{r}\|_{BMO}\leq\sqrt{1/2}\|u\|_{B}\,\sqrt{|\log(1-r^{2})|} \qquad (0<r<1),
$$
where $\|u\|_{B}=\sup_{z\in\mathbb{D}}\{|\nabla u(z)|(1-|z|^{2})\}.$
\end{Thm}

Let $BMO_{h}$ be the complex Banach space of complex-valued and
$2\pi-$periodic functions $\psi\in L^{2}(0,2\pi)$ modulo constants
with norm
$$\|\psi\|_{BMO_{h}}=\sup_{z\in\mathbb{D}}\left\{\frac{1}{2\pi}\int_{0}^{2\pi}|\psi(e^{i\theta})-f_{\psi}(z)|^{2}P(e^{i\theta},z)
\,d\theta\right\}^{1/2}<\infty.
$$
Here $P(e^{i\theta},z)$ is the Poisson kernel given by
\be\label{eq-x7}
P(e^{i\theta},z)=\frac{1-|z|^{2}}{|e^{i\theta}-z|^{2}}=\frac{\partial}{\partial
n_{\zeta}} \log\left|\frac{\zeta-z}{1-\overline{z}\zeta}\right|, \ee
where $n_{\zeta}$ is the outward normal to $\partial\mathbb{D}$ at
$\zeta=e^{i\theta}$, and
$$f_{\psi}(z)=\frac{1}{2\pi}\int_{0}^{2\pi}P(e^{i\theta},z)\psi(\theta)\,d\theta.
$$
We have several extensions of Theorem \Ref{Thm-X}.

\begin{thm}\label{thm5}
Let $\omega$ be a majorant and  $f$ be a  harmonic mapping in
$\mathbb{D}$. If $\Lambda_{f}(z)\leq
M\omega\Big(\frac{1}{1-|z|}\Big)$ in $\mathbb{D}$ and
$\psi_{r}(\theta)=f(re^{i\theta})$, then
$$\|\psi_{r}\|_{BMO}\leq2\sqrt{\omega(1)}Mr\sqrt{\int_{0}^{1}\omega\Big(\frac{1}{1-rt}\Big)\,dt},
$$
where $\theta\in[0,2\pi)$ and $r\in(0,1)$.
\end{thm}

\begin{rem}
Let $f=u+iv$ be a complex-valued continuously differentiable
function defined on $\mathbb{D}$. Then for $z=x+iy\in\mathbb{D}$,
\be\label{eqs1} \Lambda_{f}(z)\leq |\nabla u(x,y)|+|\nabla v(x,y)|,
\ee where $\nabla u=(u_{x},u_{y})$ and $\nabla v=(v_{x},v_{y})$ (see
\cite[Lemma 2]{CPW6}). Therefore, the condition $\Lambda_{f}(z)\leq
M\omega\Big(\frac{1}{1-|z|}\Big)$ in Theorem \ref{thm5} is weaker
than the condition
$$|\nabla u(x,y)|+|\nabla v(x,y)|\leq M\omega\Big(\frac{1}{1-|z|}\Big).
$$
\end{rem}

By taking $\omega(t)=t$ in Theorem \ref{thm5}, we get the following
result which  is also a generalization of \cite[Theorem 1]{KO}.

\begin{cor}\label{cor-1}
Let $f$ be a harmonic mapping in $\mathbb{D}$. If $\psi_{r}(\theta)=f(re^{i\theta})$
and $\Lambda_{f}(z)\leq M/(1-|z|)$ in $\mathbb{D}$, then
$$\|\psi_{r}\|_{BMO_{h}}\leq2\sqrt{r}M\sqrt{|\log(1-r)|},$$ where $\theta\in[0,2\pi)$ and
$r\in(0,1)$.
\end{cor}

 For $K$-quasiregular harmonic mappings with finite area, we have

\begin{cor}\label{cor-x1}
Let $f$ be a $K$-quasiregular harmonic mapping in $\mathbb{D}$
satisfying  $C=\sup_{0<r<1} S_{f}(r)<\infty$. If
$\psi_{r}(\theta)=f(re^{i\theta})$, then
$$\|\psi_{r}\|_{BMO_{h}}\leq2\sqrt{r}\sqrt{KC}\sqrt{|\log(1-r)|},$$ where $\theta\in[0,2\pi)$ and
$r\in(0,1)$.
\end{cor}

A sense-preserving and univalent harmonic mapping $f$ in
$\mathbb{D}$ will be called a {\it fully convex harmonic mapping} if
it maps every circle $|z|=r<1$ onto a convex curve  (see \cite[p.~138]{CDO}). Clunie
and Sheil-Small proved the following result.

\begin{Thm}{\rm (\cite[Corollary 5.8 ]{Clunie-Small-84})}\label{ThmA}
Let $f=h+\overline{g}$ be an univalent and sense-preserving harmonic mapping
in $\mathbb{D}$, where $h$ and $g$ are analytic in $\mathbb{D}$. If $f(\mathbb{D})$
is a convex domain, then for all $z_{1},~z_{2}\in\mathbb{D}$ with $z_{1}\neq z_{2}$,
$$|g(z_{1})-g(z_{2})|<|h(z_{1})-h(z_{2})|.$$
\end{Thm}

 In \cite{CH}, Chuaqui and
Hern\'andez discussed the relationship between the images of the
linear connectivity  under  harmonic mappings $f=h+\overline{g}$ and
under their corresponding analytic counterparts $h$, where $h$ and
$g$ are analytic in $\mathbb{D}$. For the extensive discussions on
this topic, see \cite{CPW8}. The following result is an analogous
result of Theorem \Ref{ThmA}.


\begin{thm}\label{thm6}
 Let $f=h+\overline{g}\in{\mathcal
H}$ be a fully convex harmonic mapping, where $h$ and $g$ are
analytic in $\mathbb{D}$. Then for all $r\in(0,1)$ and
$z_{1},~z_{2}\in\mathbb{D}_{r}$,
\be\label{eq-x2}
\frac{|f(z_{2})-f(z_{1})|}{1+r}\leq|h(z_{2})-h(z_{1})|\leq\frac{|f(z_{2})-f(z_{1})|}{1-r}.
\ee
\end{thm}

 The following result, which is an improvement of Theorem
\Ref{ThmA} in the case of fully convex functions,  easily follows
from Theorem \ref{thm6}.

\begin{cor}\label{cor-2}
Let $f=h+\overline{g}\in{\mathcal
H}$ be a fully convex harmonic mapping, where $h$ and $g$ are
analytic in $\mathbb{D}$. Then  $h$ is univalent in $\mathbb{D}$.
\end{cor}

We still have the following theorem due to Clunie and Sheil-Small
\cite[Theorem 5.17]{Clunie-Small-84} which helps to construct
univalent close-to-convex harmonic functions.

\begin{Thm}\label{ThmA-1}
Let $f=h +\overline{g}$ be a sense-preserving harmonic mapping in the
unit disk $\ID$, and suppose that  $h + \epsilon g$ is convex for some
$|\epsilon| \leq  1$. Then $f$ is a univalent harmonic mapping from $\ID$
onto a close-to-convex domain.
\end{Thm}

In particular, Theorem \Ref{ThmA-1} shows that a sense-preserving
harmonic mapping in $\ID$ is necessarily close-to-convex in $\ID$
whenever the analytic part of it is convex. In contrast, under a
mild restriction on $f$, namely $f_{\overline{z}}(0)=0$, Corollary
\ref{cor-2} shows that the analytic part $h$ of a sense-preserving
fully convex harmonic mapping $f=h +\overline{g}$ is necessarily
univalent in the unit disk $\ID$. On the other hand, another result
of Clunie and Sheil-Small (see for example, \cite[Theorem 5.7 and
Corollary 5.14]{Clunie-Small-84}) shows that conclusion of Corollary
\ref{cor-2}  could be improved (see \cite[Theorem 3.1]{Po}).

\section{Coefficients estimates and the Landau-Bloch theorem for locally univalent harmonic mappings}\label{csw-sec3}

\begin{Lem}{\rm (\cite[Lemma 1]{CPW0} $\mbox{or}$ \cite[Theorem 1.1]{CPW-BMMSC2011})}\label{LemgA}
Let $f$ be a harmonic mapping of $\mathbb{D}$ into $\mathbb{C}$ such
that $|f(z)|\leq M$ and
$f(z)=\sum_{n=0}^{\infty}a_{n}z^{n}+\sum_{n=1}^{\infty}\overline{b}_{n}\overline{z}^{n}$.
Then $|a_{0}|\leq M$ and for all $n\geq 1,$
$$|a_{n}|+|b_{n}|\leq \frac{4M}{\pi}.
$$
\end{Lem}

\subsection*{Proof of Theorem \ref{thm-1}} Let   $f\in {\mathcal H}(C)$.
 For $z\in\mathbb{D}$,
consider $w(z)=\overline{f_{\overline{z}}(z)}/f_{z}(z)$. Then
$w(0)=0$ and by Schwarz's lemma one has  $|w(z)|<r$ for $|z|<r$,
where $r\in(0,1)$. This implies
\be\label{eq-v}
\frac{\Lambda_{f}(z)}{\lambda_{f}(z)}=\frac{1+|w(z)|}{1-|w(z)|}\leq\frac{1+r}{1-r}:=K(r).
\ee Since $J_f(z)=\Lambda_{f}(z)\lambda_{f}(z)$ for the
sense-preserving mapping $f$, we easily have
\begin{eqnarray*}
S_{f}(r)=\int_{\mathbb{D}_{r}}\Lambda_{f}(z)\lambda_{f}(z)\,dA(z)
\geq\frac{1}{K(r)}\int_{\mathbb{D}_{r}}\Lambda_{f}^{2}(z)\,dA(z),
\end{eqnarray*}
which implies
$$\int_{\mathbb{D}_{r}}\Lambda_{f}^{2}(z)\,dA(z)\leq CK(r).
$$
For $\theta\in[0,2\pi)$ and $z\in\mathbb{D}_{r}$, let
$H_{\theta}(z)=\big(f_{z}(z)+e^{i\theta}\overline{f_{\overline{z}}(z)}\big)^{2}$.
Because $|H_{\theta}(z)|$ is subharmonic for $z\in\mathbb{D}_{r}$,
we have
\begin{eqnarray*}
|H_{\theta}(z)|&\leq&\frac{\int_{0}^{r-|z|}\int_{0}^{2\pi}|H_{\theta}(z+\rho
e^{i\beta})|\rho \,d\beta \,d\rho}{\pi(r-|z|)^{2}}\\
&\leq&\frac{1}{(r-|z|)^{2}}\int_{\mathbb{D}_{r-|z|}}(|f_{z}(z)|+|f_{\overline{z}}(z)|)^{2}\,dA(z)\\
&\leq&\frac{CK(r)}{(r-|z|)^{2}},
\end{eqnarray*}
and  the arbitrariness of $\theta\in[0,2\pi)$ gives the inequality
\be\label{thm-2e} \Lambda_{f}^{2}(z)\leq\frac{CK(r)}{(r-|z|)^{2}}.
\ee
For $\zeta\in\mathbb{D}$, let $F(\zeta)=r^{-1}f(r\zeta) .$ Then
$F(0)=0$ and by (\ref{thm-2e}) we see that \be\label{eq-crp-1}
\Lambda_{F}(\zeta)=\Lambda_{f}(z)\leq\frac{\sqrt{CK(r)}}{r}\frac{1}{1-|\zeta|},
\ee
where $z=r\zeta$. By (\ref{eq-v}) and (\ref{eq-crp-1}), we
obtain that
$$\Lambda_{F}(\zeta)\leq\min_{0<r<1}\left\{\frac{C(1+r)}{r^{2}(1-r)}\right\}^{1/2}\frac{1}{1-|\zeta|}
=\frac{Q(r_{0})}{1-|\zeta|},
$$
where
$$Q(r_{0})=\left[\frac{C(1+r_{0})}{r_{0}^{2}(1-r_{0})}\right]^{1/2}
=\sqrt{(11+5\sqrt{5})C/2} ~\approx3.330\sqrt{C}
$$
and $r_{0}=(\sqrt{5}-1)/2  \approx0.618$. Again, for
$w\in\mathbb{D}$ and a fixed $t\in(0,1)$, let $G(w)=t^{-1}F(tw)$.
Then $G(0)=0$ and
\be\label{eq-v2}
\Lambda_{G}(w)=\Lambda_{F}(\zeta)\leq\frac{Q(r_{0})}{1-|\zeta|}=\frac{Q(r_{0})}{1-t|w|}\leq\frac{Q(r_{0})}{1-t}:=M(t)
\ee
and $G$ has the form
$$G(w)=\sum_{n=1}^{\infty}A_{n}w^{n}+\sum_{n=1}^{\infty}\overline{B}_{n}\overline{w}^{n},$$
where $A_{n}=a_{n}r_{0}^{n-1}t^{n-1}$,
$B_{n}=b_{n}r_{0}^{n-1}t^{n-1}$ and $\zeta=wt$. For $\theta\in[0,2\pi)$, we consider the
function $T$ defined by
$$T(w)=G_{w}(w)+e^{i\theta}G_{\overline{w}}(w).
$$
Then $T(w)$ can be written in power series as
$$T(w)=\sum_{n=1}^{\infty}n(A_{n}w^{n-1}+e^{i\theta}\overline{B}_{n}\overline{w}^{n-1}).
$$
Applying (\ref{eq-v2}), for $w\in\mathbb{D}$, we get
$$|T(w)|\leq\Lambda_{G}(w)<M(t).$$
By  Lemma \Ref{LemgA}, for $n\geq2$, 
we have
$$n(|A_{n}+e^{i\theta}\overline{B}_{n}|)\leq \frac{4M(t)}{\pi},
$$
which gives
\begin{eqnarray*}
|a_{n}|+|b_{n}|&\leq&\frac{4Q(r_{0})}{\pi
nr_{0}^{n-1}}\inf_{0<t<1}\left\{\frac{1}{t^{n-1}(1-t)}\right\}\\
&=&\frac{4Q(r_{0})}{\pi
r_{0}^{n-1}}\left(1+\frac{1}{n-1}\right)^{n-1}\\
&<&\frac{4Q(r_{0}){\rm e}}{ \pi r_{0}^{n-1}},
\end{eqnarray*}
where ${\rm e}$ is defined by \eqref{defn-e}.
Finally, we come to prove that
$|a_{1}|\leq\sqrt{2C}$. Without loss of generality, we assume that
$$C=\int\int_{\mathbb{D}}J_{f}(z)dA(z).
$$
Then, since the dilatation $\omega $ of $f=h+\overline{g}$ satisfies
the relation $g'(z)=\omega(z) h'(z)$, by the definition of the
Jacobian, we have
\begin{eqnarray*}
C&=&\int\int_{\mathbb{D}}J_{f}(z)dA(z)\\
&=&\int\int_{\mathbb{D}}(1-|\omega(z)|^{2})|h'(z)|^{2}dA(z)\\
&\geq&\int\int_{\mathbb{D}}(1-|z|^{2})|h'(z)|^{2}dA(z) ~\mbox{ (by Schwarz' lemma $|\omega(z)|\leq |z|$)}\\
&=&\sum_{n=1}^{\infty}\frac{n}{n+1}|a_{n}|^{2}\\
&\geq&\frac{|a_{1}|^{2}}{2},
\end{eqnarray*}
which shows that $|a_{1}|\leq\sqrt{2C}$. Especially, if $C=1/2$,
then the estimate of $|a_{1}|$ is sharp and the extreme function is
$f(z)=z+\overline{z}^{2}/2$. The proof of the theorem is complete.
\qed

\subsection*{Proof of Theorem \ref{thm-1.x}}
By the hypotheses, $f=h+\overline{g}$ is a $K$-quasiregular harmonic
mapping in $\mathbb{D}$. Therefore, as in the proof of Theorem
\ref{thm-1}, we see that
$$\int_{\mathbb{D}}\Lambda_{f}^{2}(z)\,dA(z)\leq CK
$$
and
$$\Lambda_{f}^{2}(z)\leq\frac{CK}{(1-|z|)^{2}} ~\mbox{ for $z\in \ID$},
$$
so that 
\be\label{eq-e1}\Lambda_{f}(z)\leq\frac{\sqrt{CK}}{1-|z|}~\mbox{ for
$z\in \ID$}. \ee Thus, $f\in\mathcal{HB}$.

 On the other hand, for $\zeta\in\mathbb{D}$, let $F(\zeta)=r^{-1}f(r\zeta).$ For
$w\in\mathbb{D}$ and $\theta\in[0,2\pi)$, let
$$T(\zeta)=F_{\zeta}(\zeta)+e^{i\theta}F_{\overline{\zeta}}(\zeta).$$
Then
$$T(\zeta)=\sum_{n=1}^{\infty}n(a_{n}\zeta^{n-1}+e^{i\theta}\overline{b}_{n}\overline{\zeta}^{n-1})r^{n-1}
~\mbox{ for $\zeta\in \ID$},
$$
and $|T(\zeta)|<\sqrt{CK}/(1-r)$  for $\zeta\in \ID$.
 We see that
$$|T(0)|=|a_{1}+e^{i\theta}\overline{b}_{1}|<\frac{\sqrt{CK}}{1-r},$$
which implies that
$$|a_{1}|+|b_{1}|\leq\sqrt{CK}\inf_{0<r<1}\left\{\frac{1}{1-r}\right\}=\sqrt{CK}.$$
By  Lemma \Ref{LemgA}, for $n\geq2$, 
we easily have
\be\label{eq-thm-1.x}
n(|a_{n}|+|b_{n}|)\leq\frac{4\sqrt{CK}}{\pi}\inf_{0<r<1}\left\{\frac{1}{r^{n-1}(1-r)}\right\}.
\ee
Then by (\ref{eq-thm-1.x}), we get
$$|a_{n}|+|b_{n}|\leq
\begin{cases}
 \sqrt{CK} & \mbox{ if } n=1,\\[3mm]
\displaystyle
\frac{4\sqrt{CK}}{\pi}\left(1+\frac{1}{n-1}\right)^{n-1} &\mbox{ if } n\geq2.
\end{cases}
$$
The proof of the theorem is complete.
\qed

\vspace{8pt}

The following result is well-known (see for example \cite{Nehari})

\begin{Lem}\label{Lem-A}
Let $\psi$ be an analytic function in $\ID$ with $\psi (z)=\sum_{n=0}^{\infty}c_{n}z^{n}$.
If $|\psi(z)|\leq 1$, then for each $n\geq 1$, 
$|c_0|^2+|c_n|\leq 1.
$
\end{Lem}

\subsection*{Proof of Theorem \ref{thm-2}}
Let $f=h+\overline{g}\in {\mathcal H}_{\alpha}(C)$ with $0<\alpha< Q(r_{0})$,
where $r_{0}$ and $Q(r_{0})$ are the same as in Theorem \ref{thm-1}.

Following the proof of Theorem \ref{thm-1}, for $w\in\mathbb{D}$ and $\theta\in[0,2\pi)$, we let
$$H(w)=\frac{G_{w}(w)+e^{i\theta}\overline{G_{\overline{w}}(w)}}{M(t)},
$$
where $G$ and $M(t)$ are the same as in the proof of Theorem \ref{thm-1}. Then
$$H(w)=\frac{1}{M(t)}\sum_{n=1}^{\infty}n(A_{n}+e^{i\theta}B_{n})w^{n-1}
$$
and $|H(w)|<1$ for $w\in\mathbb{D}$, where $A_{n}=a_{n}r_0^{n-1}t^{n-1}$ and $B_{n}=b_{n}r_0^{n-1}t^{n-1}$.
By Lemma \Ref{Lem-A}, we have
$$\frac{n|A_{n}+e^{i\theta}B_{n}|}{M(t)}\leq1-\frac{\lambda^{2}_{G}(0)}{M^{2}(t)}=1-\frac{\alpha^{2}}{M^{2}(t)}.$$
Since $|A_{n}|+|B_{n}| =(|a_{n}| +|b_{n}|)r_0^{n-1}t^{n-1}$, the arbitrariness of $\theta\in[0,2\pi)$ gives
$$
(|a_{n}| +|b_{n}|)r_0^{n-1}t^{n-1}
\leq\frac{1}{n}\left(M(t)-\frac{\alpha^{2}}{M(t)}\right) =\frac{1}{n}\left (\frac{Q^{2}(r_{0})-\alpha^{2}(1-t)^{2}}{(1-t)Q(r_{0})}\right )
$$
which implies that
\begin{eqnarray*}
|a_{n}|+|b_{n}|&\leq& \frac{1}{nr_0^{n-1}t^{n-1}}\left (\frac{Q^{2}(r_{0})-\alpha^{2}(1-t)^{2}}{(1-t)Q(r_{0})}\right )
\\
&\leq &\frac{1}{nr_{0}^{n-1}Q(r_{0})}\inf_{0<t<1}\left\{\frac{Q^{2}(r_{0})-\alpha^{2}(1-t)^{2}}{t^{n-1}(1-t)}\right\}\\
&<&\frac{Q(r_{0})}{
nr_{0}^{n-1}}\inf_{0<t<1}\left\{\frac{1}{t^{n-1}(1-t)}\right\}
=\frac{Q(r_{0})}{
r_{0}^{n-1}}\left(1+\frac{1}{n-1}\right)^{n-1}\\
&<&\frac{Q(r_{0}){\rm e}}{ r_{0}^{n-1}},
\end{eqnarray*}
where ${\rm e}$ is defined by \eqref{defn-e}.
The proof of the theorem is complete. \qed

\subsection*{Proof of Theorem \ref{thm2.0}}
As in Theorem \ref{thm-2}, let $f=h+\overline{g}$, where $g$ and $h$
are analytic in $\mathbb{D}$ and have the form (\ref{CRP-eq-3}). For
$\zeta\in\mathbb{D}$, let $F(\zeta)=f(r_{0}\zeta)/r_{0}$, where
$r_{0}$ is the same as in Theorem \ref{thm-1}. From the proof of
Theorem \ref{thm-2},  for $n\geq2$,
we have
\be\label{eq-1-thm3} |a_{n}|+|b_{n}|
<\frac{Q(r_{0}){\rm e}}{ r_{0}^{n-1}}.
\ee
To prove the univalence of $F$, we choose two distinct points $\zeta_{1}$,
$\zeta_{2}\in\mathbb{D}_{\rho}$, where \be\label{eq-2-thm3} \rho
=1-\frac{1}{\sqrt{1+\frac{\alpha}{{\rm e} Q(r_{0})}}}. \ee
Then (\ref{eq-1-thm3}) and (\ref{eq-2-thm3}) yield that
\begin{eqnarray*}
|F(\zeta_{2})-F(\zeta_{1})|&=& \left
|\int_{[\zeta_{1},\zeta_{2}]}F_{\zeta}(\zeta)\,d\zeta+F_{\overline{\zeta}}
(\zeta)\,d\overline{\zeta}\right |\\
&\geq& \left |\int_{[\zeta_{1},\zeta_{2}]}F_{\zeta}(0)\,d\zeta+F_{\overline{\zeta}}(0)\,d\overline{\zeta}\right |\\
&& -\left
|\int_{[\zeta_{1},\zeta_{2}]}(F_{\zeta}(\zeta)-F_{\zeta}(0))\,d\zeta+(F_{\overline{\zeta}}(\zeta)-
F_{\overline{\zeta}}(0))\,d\overline{\zeta}\right |\\
&>&|\zeta_{1}-\zeta_{2}|\big[\lambda_{F}(0)-\sum_{n=2}^{\infty}(|a_{n}|+|b_{n}|)nr_{0}^{n-1}\rho^{n-1}\big]\\
&\geq& |\zeta_{1}-\zeta_{2}|\left [\alpha-Q(r_{0})\left(1+\frac{1}{n-1}\right)^{n-1}\cdot\frac{\rho(2-\rho)}{(1-\rho)^{2}}\right ]\\
&> &|\zeta_{1}-\zeta_{2}|\left [\alpha-{\rm e} Q(r_{0})
\cdot\frac{\rho(2-\rho)}{(1-\rho)^{2}}\right ]= 0, ~\mbox{ by \eqref{eq-2-thm3}}.
\end{eqnarray*}
Thus, $F(\zeta_{2})\neq F(\zeta_{1}).$ The univalence of $F$ follows
from the arbitrariness of $\zeta_1$ and $\zeta_2$. This implies that
$f$ is univalent in $\mathbb{D}_{r_{0}\rho}.$

Now, for any $\zeta'=\rho e^{i\theta}\in\partial\mathbb{D}_{\rho}$,
we easily obtain that
\begin{eqnarray*}
|F(\zeta')|&\geq&\alpha\rho-\sum_{n=2}^{\infty}(|a_{n}|+|b_{n}|)r_{0}^{n-1}\rho^{n}\\
&\geq&\alpha\rho-{\rm e}
Q(r_{0})\sum_{n=2}^{\infty}\rho^{n}\\
&=&\rho\left(\alpha-\frac{{\rm e}
Q(r_{0})\rho}{1-\rho}\right) =\frac{R_0}{r_{0}}.
\end{eqnarray*}
Therefore, $f(\mathbb{D}_{r_{0}\rho})$ contains a univalent disk of radius
$R_0$. The proof of the theorem is complete. \qed

\section{Lipschitz-type spaces of harmonic mappings}\label{csw-sec4}

\begin{Lem}$(${\rm \cite[Lemma 1]{CPW5}}$)$\label{lemma2.1}
Let $f$ be a  $K$-quasiregular harmonic mapping in $\mathbb{D}$ with
$f(\mathbb{D})\subset\mathbb{D}$. Then for all $z\in\mathbb{D},$
\be\label{eq-2.1} \Lambda_{f}(z)\leq
K\frac{1-|f(z)|^{2}}{1-|z|^{2}}. \ee Moreover, \eqref{eq-2.1} is
sharp when $K=1.$
\end{Lem}

\subsection*{Proof of Theorem \ref{thm-7}}
The implications $f\in L_{\omega}(\mathbb{D}_{r})\Rightarrow |f|\in
L_{\omega}(\mathbb{D}_{r})\Rightarrow |f|\in
L_{\omega}(\mathbb{D}_{r},\partial \mathbb{D}_{r})$ are obvious. We
only need to prove $|f|\in L_{\omega}(\mathbb{D}_{r},\partial
\mathbb{D}_{r})\Rightarrow f\in L_{\omega}(\mathbb{D}_{r}).$ For
$z\in\mathbb{D}$, let
$w(z)=\overline{f_{\overline{z}}(z)}/f_{z}(z)$. Then $w(0)=0$ and
for  any fixed $r\in(0,1)$ and  $z\in\mathbb{D}_{r}$, $|w(z)|<r$.
This gives
$$\frac{\Lambda_{f}(z)}{\lambda_{f}(z)}=\frac{1+|w(z)|}{1-|w(z)|}\leq\frac{1+r}{1-r}=K(r).$$
Now, for a fixed point $z\in \mathbb{D}_{r}$, we consider the function
$$F(\eta)=f(z+d(z)\eta)/M_{z},\ \eta\in\mathbb{D},
$$
where $d(z):=d(z,\partial
\mathbb{D}_{r})$ denotes the Euclidean distance from $z$ to the
boundary $\partial\mathbb{D}_{r}$ of $\mathbb{D}_{r}$  and
$M_{z}:=\sup\{|f(\zeta)|:\, |\zeta-z|<d(z)\}.$ By an elementary
calculation, we obtain that
$$\frac{\Lambda_{F}(\eta)}{\lambda_{F}(\eta)}=\frac{\Lambda_{f}(\xi)}{\lambda_{f}(\xi)}
\leq K(r),
$$
where $\xi=z+d(z)\eta.$ Then $F$ is a  $K(r)$-quasiregular harmonic
mapping of $\mathbb{D}$ into itself. By Lemma \Ref{lemma2.1}, we see that
$$\Lambda_{F}(0)\leq K(r)(1-|F(0)|^{2})
$$
which may be written as \be\label{eqh4}
d(z)\Lambda_{f}(z)\leq2K(r)(M_{z}-|f(z)|). \ee Without loss of
generality, we let $\zeta\in\partial \mathbb{D}_{r}$ with
$|\zeta-z|=d(z)$, and let $w\in\mathbb{D}(z,d(z))$. Then
\begin{eqnarray*}
|f(w)|-|f(z)|&\leq&\big ||f(w)|-|f(\zeta)|\big |+\big ||f(\zeta)|-|f(z)|\big |\\
&\leq&M\omega(d(z))+M\omega(2d(z))\\
&\leq&3M\omega(d(z)),
\end{eqnarray*}
and thus
$$\sup_{w\in\mathbb{D}(z,d(z))}(|f(w)|-|f(z)|)\leq 3M\omega(d(z)),
$$
which implies that
$M_{z}-|f(z)|\leq 3M\omega(d(z)).$
This inequality together with (\ref{eqh4}) shows that
\be\label{eq4.00x}
\Lambda_{f}(z)\leq6MK(r)\cdot\frac{\omega(d(z))}{d(z)},\ z\in
\mathbb{D}_{r}.
\ee

Finally, given any two points $z_{1},z_{2}\in \mathbb{D}_{r}$, let
$\gamma\subset \mathbb{D}_{r}$ be a rectifiable  curve  joining
$z_{1}~\mbox{and}~z_{2}$. Integrating (\ref{eq4.00x}) along
$\gamma$, we obtain that \be\label{eq6h}
|f(z_{1})-f(z_{2})|\leq\int_{\gamma}(|f_{z}(z)|+|f_{\overline{z}}(z)|)\,ds(z)
\leq6MK(r)\int_{\gamma}\frac{\omega(d(z))}{d(z)}\,ds(z). \ee
Therefore, (\ref{eq6h}) yields
$$|f(z_{1})-f(z_{2})|\leq M_{1}\cdot\omega(|z_{1}-z_{2}|),
$$
where $M_{1}$ is a positive constant. The proof of the theorem is
complete. \qed


\subsection*{Proof of Theorem \ref{thm5}}
By (\ref{eq-x7}) and integrating by parts, we have

\vspace{8pt}

$\ds \frac{1}{2\pi}\int_{0}^{2\pi}|f(re^{i\theta})-f(rz)|^{2}P(e^{i\theta},z)\,d\theta
$

\begin{eqnarray*}
&=&\frac{1}{2\pi}\int_{0}^{2\pi}|f(re^{i\theta})-f(rz)|^{2}\frac{\partial}{\partial
n_{\zeta}}\log\left|\frac{\zeta-z}{1-\overline{z}\zeta}\right|\,d\theta\\
&=&-\frac{2r^{2}}{\pi}\int_{\mathbb{D}}\log\left|\frac{w-z}{1-\overline{z}w}\right|
\left(|f_{\xi}(rw)|^{2}+|f_{\overline{\xi}}(rw)|^{2}\right)\,dA(w)\\
&\leq&-\frac{2r^{2}}{\pi}\int_{\mathbb{D}}\log\left|\frac{w-z}{1-\overline{z}w}\right|
\Lambda_{f}^{2}(rw)\,dA(w)\\
&\leq&-
\frac{2r^{2}M^{2}}{\pi}\int_{\mathbb{D}}\log\left|\frac{w-z}{1-\overline{z}w}\right|\omega^{2}\Big(\frac{1}{1-r|w|}\Big)
\,dA(w)\\
&\leq&\frac{2r^{2}M^{2}}{\pi}\int_{\mathbb{D}}\log\frac{1}{|w|}\omega^{2}\Big(\frac{1}{1-r|w|}\Big)
\,dA(w)\\
&=&4r^{2}M^{2}\int_{0}^{1}\rho|\log\rho|\omega^{2}\Big (\frac{1}{1-r\rho}\Big)\,d\rho\\
&=&4r^{2}M^{2}\int_{0}^{1}\left[|\log\rho|\frac{d}{d\rho}\int_{0}^{\rho}t\omega^{2}\Big(\frac{1}{1-rt}\Big)\,dt\right]\,d\rho\\
&=&4r^{2}M^{2}\left[|\log\rho|\int_{0}^{\rho}t\omega^{2}\Big(\frac{1}{1-rt}\Big)\,dt\Big|_{0}^{1}
+\int_{0}^{1}\frac{\int_{0}^{\rho}t\omega^{2}\Big(\frac{1}{1-rt}\Big)\,dt}{\rho}\,d\rho\right]\\
&=&4r^{2}M^{2}\int_{0}^{1}\frac{\int_{0}^{\rho}t\omega^{2}\Big(\frac{1}{1-rt}\Big)\,dt}{\rho}\,d\rho\\
&\leq&4r^{2}M^{2}\int_{0}^{1}\int_{0}^{\rho}\omega^{2}\Big(\frac{1}{1-rt}\Big)\,dt\,d\rho\\
&=&4r^{2}M^{2}\int_{0}^{1}\left(\int_{t}^{1}d\rho\right)\omega^{2}\Big(\frac{1}{1-rt}\Big)\,dt\\
&=&4r^{2}M^{2}\int_{0}^{1}(1-t)\omega^{2}\Big(\frac{1}{1-rt}\Big)\,dt\\
&=&4r^{2}M^{2}\int_{0}^{1}\frac{1-t}{1-rt}\left[(1-rt)\omega\Big(\frac{1}{1-rt}\Big)\right]\omega\Big(\frac{1}{1-rt}\Big)\,dt\\
&\leq&4r^{2}M^{2}\omega(1)\int_{0}^{1}\omega\Big(\frac{1}{1-rt}\Big)\,dt,
\end{eqnarray*} which implies
$$\|\psi_{r}\|_{BMO_{h}}\leq2\sqrt{\omega(1)}Mr\sqrt{\int_{0}^{1}\omega\Big(\frac{1}{1-rt}\Big)\,dt},
$$
where $\xi=rw$. The proof  is complete.  \qed

 \subsection*{Proof of Corollary \ref{cor-x1}} Corollary  \ref{cor-x1} easily follows from
(\ref{eq-e1}) and Corollary  \ref{cor-1}. \qed

\subsection*{Proof of Theorem \ref{thm6}}  Differentiating both sides of equation $f^{-1}(f(z))=z$ yields
the relations
$$(f^{-1})_{\zeta}h'+(f^{-1})_{\overline{\zeta}}g'=1
~\mbox{ and } ~(f^{-1})_{\zeta}\overline{g'}+(f^{-1})_{\overline{\zeta}}\overline{h'}=0,
$$
where $\zeta=f(z)$. This gives \be\label{eq-x1}
 (f^{-1})_{\zeta}=\frac{\overline{h'}}{J_{f}}~\mbox{and}~(f^{-1})_{\overline{\zeta}}
=-\frac{\overline{g'}}{J_{f}}.
\ee
Since $\Omega=f(\mathbb{D}_{r})$ is  convex, for any two distinct points
$z_{1},~z_{2}\in\mathbb{D}_{r}$ and $t\in[0,1]$, we have
$$\varphi(t)=(f(z_{2})-f(z_{1}))t+f(z_{1})\in\Omega.
$$
Let $\gamma=f^{-1}\circ\varphi$ and
$f(z_{2})-f(z_{1})=|f(z_{2})-f(z_{1})|e^{i\theta_{0}}$. For
$z\in\mathbb{D}$, let
$w(z)=\overline{f_{\overline{z}}(z)}/f_{z}(z)$. Then $w(0)=0$ and
for $z\in\mathbb{D}_{r}$, $|w(z)|<r$, where $r\in(0,1)$. This
implies that $f$ is $K(r)$-quasiconformal harmonic mapping in
$\mathbb{D}_{r}$, where $K(r)=\frac{1+r}{1-r}$.
By calculations and (\ref{eq-x1}), we have

\vspace{6pt}
\noindent

$|h(z_{2})-h(z_{1})|$
\begin{eqnarray*}
&=&\left|\int_{\gamma}h'(z)\,dz\right|\\
&=&\left|\int_{0}^{1}h'(\gamma(t))\,\frac{d}{dt}\gamma(t)\,dt\right|\\
&=&\left|\int_{0}^{1}h'(\gamma(t))\Big[\varphi'(t)\frac{\partial}{\partial
\zeta}f^{-1}(\varphi(t))+\overline{\varphi'(t)}\frac{\partial}{\partial
\overline{\zeta}}f^{-1}(\varphi(t))\Big]\,dt\right|\\
&=&|f(z_{2})-f(z_{1})|\left|\int_{0}^{1}h'(\gamma(t))\Big(\frac{\overline{h'(\gamma(t))}}{J_{f}(\gamma(t))}e^{i\theta_{0}}-
\frac{\overline{g'(\gamma(t))}}{J_{f}(\gamma(t))}e^{-i\theta_{0}}\Big)\,dt\right|\\
&\leq&|f(z_{2})-f(z_{1})|\int_{0}^{1}|h'(\gamma(t))|\frac{(|h'(\gamma(t))|+|g'(\gamma(t))|)}{J_{f}(\gamma(t))}\,dt\\
&=&|f(z_{2})-f(z_{1})|\int_{0}^{1}\frac{1}{1-|w(\gamma(t))|}\,dt\\
&\leq&\frac{|f(z_{2})-f(z_{1})|}{1-r}.
\end{eqnarray*}
This gives the second inequality in (\ref{eq-x2}). Next we prove the
first inequality. Applying (\ref{eq-x1}), we see that

\begin{eqnarray*}
&&\mbox{Re}\,[e^{-i\theta_{0}}(\overline{g(z_{2})-g(z_{1})})]
=\mbox{Re}\left[e^{-i\theta_{0}}\Big(\overline{\int_{0}^{1}g'(\gamma(t))\,\frac{d}{dt}\gamma(t)\,dt}\Big)\right]\\
&=&\mbox{Re}\left\{e^{-i\theta_{0}}\Big[\overline{\int_{0}^{1}g'(\gamma(t))\Big(\varphi'(t)\frac{\partial}{\partial
\zeta}f^{-1}(\varphi(t))+\overline{\varphi'(t)}\frac{\partial}{\partial
\overline{\zeta}}f^{-1}(\varphi(t))\Big)\,dt}\Big]\right\}\\
&=&|f(z_{2})-f(z_{1})|\mbox{Re}\left\{e^{-i\theta_{0}}
\Big[\overline{\int_{0}^{1}\frac{g'(\gamma(t))\overline{h'(\gamma(t))}e^{i\theta_{0}}-
|g'(\gamma(t))|^{2}e^{-i\theta_{0}}}
{J_{f}(\gamma(t))}\,dt}\Big]\right\}\\
&\leq&|f(z_{2})-f(z_{1})|\int_{0}^{1}\frac{|h'(\gamma(t))\overline{g'(\gamma(t))}e^{-2i\theta_{0}}|-
|g'(\gamma(t))|^{2}} {J_{f}(\gamma(t))}\,dt\\
&\leq&|f(z_{2})-f(z_{1})|\int_{0}^{1}\frac{|w(\gamma(t))|}{1+|w(\gamma(t))|}\,dt\\
&\leq&\frac{r|f(z_{2})-f(z_{1})|}{1+r},
\end{eqnarray*}
which gives
\be\label{eq-x3}
\mbox{Re}\left\{\frac{\overline{g(z_{2})-g(z_{1})}}{f(z_{2})-f(z_{1})}\right\}\leq\frac{r}{1+r}.
\ee

It is not difficult to see that \be\label{eq-x4}
\mbox{Re}\left\{\frac{h(z_{2})-h(z_{1})}{f(z_{2})-f(z_{1})}\right\}=1-
\mbox{Re}\left\{\frac{\overline{g(z_{2})-g(z_{1})}}{f(z_{2})-f(z_{1})}\right\}.
\ee By (\ref{eq-x3}) and (\ref{eq-x4}), we get

\begin{eqnarray*}
\frac{|h(z_{2})-h(z_{1})|}{|f(z_{2})-f(z_{1})|}&\geq&\mbox{Re}\left\{\frac{h(z_{2})-h(z_{1})}{f(z_{2})-f(z_{1})}\right\}\\
&=&1-
\mbox{Re}\left\{\frac{\overline{g(z_{2})-g(z_{1})}}{f(z_{2})-f(z_{1})}\right\}\\
&\geq&1-\frac{r}{1+r}\\
&=&\frac{1}{1+r}.
\end{eqnarray*}
The proof of this theorem is complete. \qed

\bigskip

{\bf Acknowledgements:} The authors would like to thank the referee
for his (or her) careful reading of this paper and useful
suggestions. This research was partly supported by the Construct
Program of the Key Discipline in Hunan Province and the Start
Project of Hengyang Normal University (No. 12B34).
The second author is on leave from the Department of Mathematics,
Indian Institute of Technology Madras, Chennai-600 036, India

\bibliographystyle{line}

\end{document}